\documentclass[10pt]{amsart}
\usepackage{amssymb,latexsym}
\usepackage{psfig}
\topskip  \numberwithin{equation}{section}
\newtheorem{theorem}{Theorem}[section]
\newtheorem{proposition}[theorem]{Proposition}
\newtheorem{corollary}[theorem]{Corollary}

\newtheorem{lemma}[theorem]{Lemma}

\newtheorem{example}[theorem]{Example}

\newcommand{\psdraw}[3]{\begin{array}{c} \hspace{-1mm}
\raisebox{-4pt}{\psfig{figure=#1.ps,width=#2,height=#3}}
\hspace{-1mm}\end{array}}

\def\pa{{k(k+2)(k+3)\cdots(2k-1)1(2k)23\cdots(k+1)}}
\def\pb{{k(k+2)(k+3)\cdots(2k-1)(2k)123\cdots(k+1)}}
\def\pc{{(k+1)(k+2)(k+3)\cdots(2k-1)1(2k)23\cdots k}}
\def\pd{{(k+1)(k+2)(k+3)\cdots(2k-1)(2k)123\cdots k}}

\def\sl{\sum\limits}

\begin{document}
\def\PP{\mathfrak P}
\def\ttx{\left(\frac{1}{2\sqrt{x}}\right)}
\def\SS{\frak S}

\pagenumbering{arabic} \pagestyle{headings}

\title{$321$-polygon-avoiding permutations and Chebyshev polynomials}

\maketitle

\begin{center}{{\sc Toufik Mansour}\\
\smallskip
              LaBRi, Universit\'e Bordeaux I,\\
              351 cours de la Lib\'eration, 33405 Talence Cedex\\
        {\tt toufik@labri.fr}\\
\bigskip
{\sc Zvezdelina Stankova}\\
\smallskip
              Mills College, Oakland, CA\\
 {\tt stankova@mills.edu}}
\end{center}
\begin{abstract}
A {\em $321$-$k$-gon-avoiding permutation} $\pi$ avoids $321$ and
the following four patterns:
$$\begin{array}{l}
\pa,\quad \pb,\\
\pc,\quad \pd.
\end{array}$$ The $321$-$4$-gon-avoiding
permutations were introduced and studied by Billey and Warrington
\cite{BW} as a class of elements of the symmetric group whose
Kazhdan-Lusztig, Poincar\'e polynomials, and the singular loci of
whose Schubert varieties have fairly simple formulas and
descriptions. Stankova and West \cite{SW} gave an exact
enumeration in terms of linear recurrences with constant
coefficients for the cases $k=2,3,4$.  In this paper, we extend
these results by finding an explicit expression for the generating
function for the number of $321$-$k$-gon-avoiding permutations on
$n$ letters. The generating function is expressed via Chebyshev
polynomials of the second kind.
\end{abstract}
\section{Introduction}
\noindent {\bf Definition 1.} Let $\alpha\in S_n$ and $\tau\in
S_k$ be two permutations. Then $\alpha$ {\it contains\/} $\tau$ if
there exists a subsequence $1\leq i_1<i_2<\dots<i_k\leq n$ such
that $(\alpha_{i_1}, \dots,\alpha_{i_k})$ is order-isomorphic to
$\tau$; in such a context $\tau$ is usually called a {\it
pattern\/}; $\alpha$ {\it avoids\/} $\tau$, or is $\tau$-{\it
avoiding\/}, if $\alpha$ does not contain such a subsequence. The
set of all $\tau$-avoiding permutations in $S_n$ is denoted by $
S_n(\tau)$. For a collection of patterns $T$, $\alpha$ {\it
avoids} $T$ if $\alpha$ avoids all $\tau\in T$; the corresponding
subset of $S_n$ is denoted by $ S_n(T)$.

While the case of permutations avoiding a single pattern has
attracted much attention, the case of multiple pattern avoidance
remains less investigated. In particular, it is natural to
consider permutations avoiding pairs of patterns $\tau_1$,
$\tau_2$. This problem was solved completely for $\tau_1,\tau_2\in
S_3$ (see \cite{SS}), for $\tau_1\in S_3$ and $\tau_2\in S_4$ (see
\cite{W}), and for $\tau_1,\tau_2\in S_4$ (see \cite{B1,Km} and
references therein). Several recent papers
\cite{CW,MV1,Kr,MV2,MV3,MV4} deal with the case $\tau_1\in S_3$,
$\tau_2\in S_k$ for various pairs $\tau_1,\tau_2$. The tools
involved in these papers include Fibonacci numbers, Catalan
numbers, Chebyshev polynomials, continued fractions, and Dyck
words, e.g. in \cite{MV2}:

\begin{theorem}[Mansour, Vainshtein] Let
$U_m(\cos\theta)=\sin(m+1)\theta/\sin\theta$ be the Chebyshev
polynomial of the second kind. When $2\leq d+1\leq k$, the
generating function for the number of permutations in $
S_n(321,(d+1)\cdots k12\cdots d)$ is given by
    $$\dfrac{U_{k-1}\ttx}{\sqrt{x}U_k\ttx}\cdot$$
\end{theorem}

\medskip
Recently, a special class of restricted permutations has arisen in
representation theory.

\medskip
\noindent{\bf Definition 2.} A permutation $\pi$ is {\em
$k$-gon-avoiding} if it avoids each pattern in the set $\PP_k$:
\[\begin{array}{l}\big\{ \pa,\quad \pb,\\\,\,\,\pc,\quad \pd \big\}.
\end{array} \]
We say that $\pi$ is a {\em $321$-$k$-gon}-avoiding permutation if
it is both $k$-gon-avoiding and $321$-avoiding. The number of
$321$-$k$-gon-avoiding permutations in $ S_n$ is denoted by
$f_k(n)$. The corresponding generating function is
$f_k(x)=\sum_{n\geq0} f_k(n)x^n$.

\begin{figure}[h]
$$\psdraw{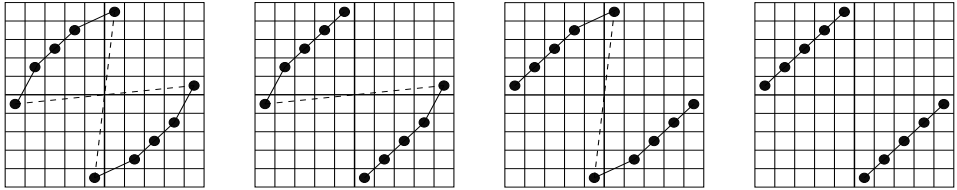}{5in}{1in}$$
\caption{$\mathfrak{B}_5$: all four 5-gons}
\end{figure}

\noindent Note that $f_k(n)=\frac{1}{n+1}\binom{2n}{n}$ for
$n\in[0,2k-1]$, as these count the permutations in $S_n(321)$ (see
\cite{Kn}).

\medskip
Billey and Warrington \cite{BW} introduced the 321-4-gon-avoiding
(or {\it 321-hexagon-avoiding}) permutations as a class in $ S_n$
whose Kazhdan-Lusztig and Poincar\'e polynomials, and the singular
loci of whose Schubert varieties have fairly simple formulas and
descriptions. Upon their request, Stankova and West \cite{SW}
presented an exact enumeration for the cases $k=2,3,4$ by using
generating trees, the symmetries in the set of the $\PP_k$, and
the structure of the $321$-avoiding permutations via Schensted's
321--subsequences decomposition.

\begin{theorem}[Stankova,West]
For $k=2,3,4$, the sequences $f_k(n)$ satisfy the recursive
relations
\[\begin{array}{l}
f_4(n)=6f_4(n-1)-11f_4(n-2)+9f_4(n-3)-4f_4(n-4)-4f_4(n-5)+f_4(n-6),\,\,n\geq 6;\\
f_3(n)=4f_3(n-1)-4f_3(n-2)+3f_3(n-3)+f_3(n-4)-f_3(n-5),\,\,n\geq 5;\\
f_2(n)=3f_2(n-1)-3f_2(n-2)+f_2(n-3)=(n-1)^2+1,\,\,n\geq 3.\\
\end{array}\]
\end{theorem}

\smallskip
In this paper we present an approach to the study of
$321$-$k$-gon-avoiding permutations in $ S_n$ which generalizes
the methods in \cite{SW} and \cite{MV3}. As a consequence, we
extend the results in \cite{SW} to all 321-$k$-gon-avoiding
permutations, and derive a number of other related results. The
main theorem of the paper is formulated as follows.

\begin{theorem}\label{thm}
For $k\geq 3$ and $s\geq 1$, define
$L_n^k(s)=\sl_{j=0}^{s}(-1)^j\binom{s-j}{j}f_k(n-j)$. When $n\geq
2k$, this sequence satisfies the linear recursive relation with
constant coefficients:
\begin{eqnarray*}
L_n^k(2k-2)=L_{n-1}^k(2k-2)+L_{n-3}^k(2k-5)+L_{n-3}^k(2k-4)+
L_{n-4}^k(2k-5)+L_{n-4}^k(2k-4).\\
\end{eqnarray*}
\end{theorem}

\begin{corollary}\label{thmm}
For $k\geq 3$,
$$f_k(x)=\dfrac{(1+x^2-x^4)U_{2k-3}\ttx-\sqrt{x}(1+x^3)U_{2k-4}\ttx}
{\sqrt{x}\left[(1+x^2-x^4)U_{2k-2}\ttx-\sqrt{x}(1+x^3)U_{2k-3}\ttx\right]}\cdot$$
\end{corollary}

The proofs of Theorem~\ref{thm} and Corollary~\ref{thmm} are
presented in Section~2. Note that Corollary~\ref{thmm} implies the
previously known results for the cases $k=3,4$ (see \cite{SW}):
\begin{eqnarray*}
f_3(x)&=&\dfrac{1-3x+2x^2-2x^3-2x^4+3x^5}{1-4x+4x^2-3x^3-x^4+5x^5-x^6},\\
f_4(x)&=&\dfrac{1-5x+7x^2-5x^3+x^4+7x^5-4x^6}{1-6x+11x^2-9x^3+4x^4+8x^5-9x^6+x^7}\cdot
\end{eqnarray*}

In Section~3, we describe several generalizations of
Theorem~\ref{thm} and Corollary~\ref{thmm}, following similar
arguments from their proofs.
\section{Proof Theorem~\ref{thm}}

\subsection{Refinement of the numbers $f_k(n)$}
$\phantom{10mm}$

\noindent{\bf Definition 3.} For $m,n$ with $1\leq m\leq n$ and
distinct $i_1,i_2,...,i_m\in\mathbb{N}$, we denote by $f_k(n;
i_1,\dots,i_m)$ the number of $321$-$k$-gon-avoiding permutations
$\pi\in S_n$ which start with $i_1i_2\cdots i_m$:
$\pi_1\pi_2\cdots\pi_m=i_1i_2\cdots i_m$. The corresponding subset
of $ S_n$ is denoted by ${\mathcal F}_k(n; i_1,\dots,i_m)$.

\smallskip
Here follow basic properties of the numbers $f_k(n;
i_1,\dots,i_m)$, easily deduced from the definitions.

\begin{lemma}\label{l22} Let $n\geq 3$, $1\leq m\leq n$
and $1\leq i_1<i_2<\cdots < i_m\leq n$.
\begin{itemize}
\item[(a)] If $m\leq n-2$ and $3\leq i_1$, then $f_k(n;
i_1,\dots,i_m,j)=0$ for $2\leq j\leq i_m-1$. Consequently,
\[\displaystyle{f_k(n; i_1,\dots,i_m)=f_k(n;
i_1,\dots,i_m,1)+\sl_{j=i_m+1}^nf_k(n; i_1,\dots,i_m,j).}\]

\item[(b)] If $m\leq k-3$ and $2\leq i_1$, then
$f_k(n; i_1,\dots,i_m,1)=f_k(n-1; i_1-1,\dots,i_m-1).$

\smallskip
\item[(c)] If $i_1\leq k-1$,
then $f_k(n; i_1,\dots,i_m)=f_k(n-1; i_2-1,\dots,i_m-1).$
\end{itemize}
\end{lemma}
\begin{proof}
For (a), observe that if $\pi\in {\mathcal F}_k(n;
i_1,\dots,i_m,j)$ then the entries $i_m$, $j$, $1$ give an
occurrence of $321$ in $\pi$.  For the second part of (a),
consider the entries $\pi_{m+1}$ of $\pi$.  Again, avoiding $321$
forces $\pi_{m+1}=1$ or $\pi_{m+1}>i_m$.

For (b), denote by $\pi'$ the permutation obtained from $\pi$ by
deleting its smallest entry and decreasing all other entries by
$1$. Then $\pi\in {\mathcal F}_k(n; i_1,\dots,i_m,1)$ if and only
if $\pi'\in {\mathcal F}_k(n-1; i_1,\dots,i_m)$, since entry $1$
placed as in (b) cannot be used in an occurrence of $321$ or
$\tau\in\PP_k$ in $\pi$.

For (c), observe that if $\pi_1\pi_2\cdots\pi_m=i_1i_2\cdots i_m$
then the entry $i_1$ cannot appear in any occurrences of
$\tau\in\PP_k$; further, if there is an occurrence $xyz$ of $321$
such that $x=i_1$ then there is an occurrence $i_2yx$ of $321$ in
$\pi$.
\end{proof}

Lemma~\ref{l22} implies an explicit formula for $f_k(n; s)$ for
the first values of $s$.

\begin{proposition}\label{por1}
For $1\leq s\leq \min\{k-1,n\}$,
    $$f_k(n; s)=\sum_{j=0}^{s-1}(-1)^j\binom{s-1-j}{j}f_k(n-1-j).$$
\end{proposition}
\noindent{\it Proof.}  By Lemma~\ref{l22}(c)-(a) the proposition
holds for $s=1,2$. For $s\geq 3$, Lemma~\ref{l22}(b) says
\begin{eqnarray*}
    f_k(n;s)&=&f_k(n;s,1)+\sum_{j=s+1}^n f_k(n;s,j)\\
    \Rightarrow\,\,f_k(n;s)&=&f_k(n-1;s-1)+\sum_{j=s+1}^{n}
f_k(n-1;j-1)=\sum_{j=s-1}^{n-1}f_k(n-1;j).
\end{eqnarray*}
Equivalently,
\begin{equation}
f_k(n;s)=f_k(n-1)-\sum_{j=1}^{s-2} f_k(n-1;j). \label{eqa1}
\end{equation}
Using induction on $s$, we assume that the proposition holds for
all $1\leq j\leq s-1$. Then (\ref{eqa1}) yields
    $$f_k(n;s)=f_k(n-1)-\sum_{j=1}^{s-2}\sum_{i=0}^{j-1}(-1)^i\binom{j-1-i}{i}f_k(n-2-i).$$
>From the familiar equality
$\binom{1}{a}+\binom{2}{a}+\cdots+\binom{b}{a}=\binom{b+1}{a+1}$
we obtain for all $2\leq s\leq k-1$
\begin{eqnarray*}
f_k(n;s)&=&f_k(n-1)+\sum_{j=1}^{s-1}(-1)^j\binom{s-1-j}{j}f_k(n-1-j)\\
&=&\sum_{j=0}^{s-1}(-1)^j\binom{s-1-j}{j}f_k(n-1-j).  \qed
\end{eqnarray*}

Next we introduce objects $A_d(n,m)$ which organize suitably the
information about the numbers $f_k(n;i_1,i_2,...,i_m)$ and play an
important role in the proof of the main result.

\smallskip
\noindent{\bf Definition 4.} For $1\leq d\leq n+1-m$ and $1\leq
m\leq n$ set
\begin{equation*}
A_d(n,m)=\sum_{d\leq i_1<i_2<\cdots < i_m\leq n} f_k(n;
i_1,\dots,i_m).\label{obja}
\end{equation*}
In the following subsections \ref{first}--\ref{second} we derive
two expressions for $A_k(n,k-1)$, compare them in subsection
\ref{third}, and thus complete the proof of Theorem \ref{thm}.

\subsection{First expression for $A_k(n,k-1)$}\label{first} The numbers
$A_d(n,m)$ satisfy the following recurrence.

\begin{lemma}\label{lem1}
For $2\leq d\leq k$ and $1\leq m\leq \min\{k-1,n\}$,
        $$A_d(n,m)=A_{d-1}(n,m)-A_{d-1}(n-1,m-1).$$
\end{lemma}
\noindent{\sc Proof.} By (\ref{obja}), for all $2\leq  d\leq  k$
we have:
$$A_d(n,m)=A_{d-1}(n,m)-\sum_{d\leq i_2<\cdots <i_m\leq
n}f_k(n;d-1,i_2,...,i_m).$$ Lemma~\ref{l22}(c) and Definition 4
imply
\begin{eqnarray*}
A_d(n,m)&=&A_{d-1}(n,m)-\sum_{d-1=i_2<i_3<\cdots < i_m\leq n-1}
f_k(n-1; i_2,i_3\dots,i_m)\\
&=& A_{d-1}(n,m)-A_{d-1}(n-1,m-1). \qed
\end{eqnarray*}

We next find an explicit expression for $A_1(n,m)$ in terms of
$f_k(n)$.

\begin{lemma}\label{lem2}
Let $1\leq  m\leq  \min\{k-1,n\}$. Then
    $$A_1(n,m)=\sum_{j=0}^m (-1)^j\binom{m-j}{j}f_k(n-j).$$
\end{lemma}
\noindent{\it Proof.} For $m=1$, $A_1(n,1)=\sum_{1\leq i_1\leq
n}f_k(n;i_1)= f_k(n)$, which equals the required expression.
Assume the lemma for $m$ and all appropriate $n$. Comparing the
$(m+1)^{\text{st}}$ entry $j$ of $\pi$ with $i_m$,
\[A_1(n,m)=A_1(n,m+1)
+\sum_{1\leq i_1<i_2<...<i_m\leq n} \sum_{j=1}^{i_m-1}
f_k(n;i_1,i_2,...,i_m,j).\] For the summation part on the
right--hand side, avoidance of (321) implies that all numbers
$1,2,...,j-1$ appear before the entry $j$, and hence $j\leq m$.
From Lemma~\ref{l22},
\begin{eqnarray*}
A_1(n,m)&=&A_1(n,m+1) +\sum_{j=1}^{m}\sum_{j+1\leq
i_j<i_{j+1}<...<i_m\leq n} f_k(n;1,2,...,j-1,i_j,i_{j+1},...,i_m,j)\\
&=&A_1(n,m+1) +\sum_{j=1}^{m}\sum_{2\leq
i_j<i_{j+1}<...<i_m\leq n-j+1} f_k(n-j+1;i_j,i_{j+1},...,i_m,1)\\
&=&A_1(n,m+1) +\sum_{j=1}^{m}\sum_{1\leq
i_j<i_{j+1}<...<i_m\leq n-j} f_k(n-j;i_j,i_{j+1},...,i_m)\\
&=&A_1(n,m+1)+\sum_{j=1}^{m} A_1(n-j,m-j+1)= \sum_{j=0}^{m}
A_1(n-j,m-j+1)
\end{eqnarray*}
Applying the above to $A_1(n-1,m-1)$, subtracting the results, and
using the induction hypothesis, we arrive at
\begin{eqnarray*}A_1(n,m)&=&A_1(n,m+1)+A_1(n-1,m-1)\\
\Rightarrow\,\,A_1(n,m+1)&=&A_1(n,m)-A_1(n-1,m-1)\\
&=&\sl_{j=0}^{m} (-1)^j\binom{m-j}{j}f_k(n-j)-
\sl_{j=0}^{m-1} (-1)^j\binom{m-j-1}{j}f_k(n-1-j)\\
&=&\sum_{j=0}^{m+1} (-1)^j\binom{m-j+1}{j}f_k(n-j). \qed
\end{eqnarray*}

Lemmas~\ref{lem1}--\ref{lem2} can be combined to derive the
following explicit expression for $A_d(n,m)$, which is easily
proven by induction on $d$.

\begin{corollary}\label{thss1}
Let $1\leq d\leq k$, $1\leq m\leq \min\{k-1,n\}$. Then
    $$A_d(n,m)=\sum_{j=0}^{d+m-1}(-1)^j\binom{d+m-1-j}{j}f_k(n-j)=
L_n^k(d+m-1).$$
\end{corollary}

In particular,
\[A_k(n,k-1)=\sum_{j=0}^{2k-2}(-1)^j\binom{2k-2-j}{j}f_k(n-j)=L_n^k(2k-2).\]

\subsection{Second expression for $A_k(n,k-1)$}\label{second}
We start by introducing three objects related to $A_d(n,m)$.

\smallskip
\noindent{\bf Definition 5.} For $1\leq d\leq n-m+1$ and $1\leq
m\leq n$ set
\begin{eqnarray*}
B_d(n,m)&=&\sum_{d+1\leq i_1<i_2<\cdots <i_m\leq n}f_k(n;
d,i_1,\dots,i_m);\\ C_d(n,m)&=& \,\,\,\sum_{d\leq i_1<i_2<\cdots
<i_m\leq n} \,\,\,f_k(n; i_1,\dots,i_m,1);\\ D_d(n,m)&=&
\sum_{d+1\leq i_1<i_2<\cdots <i_m\leq n}f_k(n; d,i_1,\dots,i_m,1).
\end{eqnarray*}
Note that by Lemma~\ref{l22}(a), for $k\geq 2$:
\begin{equation}
A_k(n,k-1)=A_k(n,k)+C_k(n,k-1).\label{eqb1}
\end{equation}
The following Propositions~\ref{prob1}--\ref{prob2} describe
$A_k(n,k)$ and $C_k(n,k-1)$ in terms of $f_k(n)$.

\begin{proposition}\label{prob1}
Let $n\geq k-1$. Then
\begin{itemize}
\item[(a)] $C_k(n,k-1)=C_{k-1}(n-1,k-1)+A_{k-1}(n-3,k-2)$;

\item[(b)] $C_k(n,k-1)=C_k(n-1,k-1)+A_{k-2}(n-3,k-2)+A_{k-1}(n-3,k-2)$.
\end{itemize}
\end{proposition}
\begin{proof}
For (a), similarly to Lemma~\ref{l22}(b) (with $k\geq3$), we have
\begin{equation}
C_k(n,k-1)=\sum_{k\leq i_1<\cdots < i_{k-1}\leq n}
\!\!\!\!\!\!f_k(n; i_1,\dots,i_{k-1},1,2)+\sum_{k\leq i_1<\cdots <
i_{k-1}<i_k\leq n} \!\!\!\!\!\!f_k(n; i_1,\dots,i_{k-1},1,i_k).
\label{eqsa}
\end{equation}
If $\pi$ starts with $i_1,\dots,i_{k-1},1,2$ as in the first sum
in (\ref{eqsa}), then the entry 2 cannot participate in an
occurrence of $321$ or of $\tau\in\PP_k$. Hence
\begin{equation}
f_k(n;i_1,\dots,i_{k-1},1,2)=f_k(n-1;i_1-1,\dots,i_{k-1}-1,1).\label{eqaa1}
\end{equation}
For the second sum in (\ref{eqsa}), if $\pi$ starts with
$i_1,\dots,i_{k-1},1,i_k$, avoidance of 321 and both
$$\pc\,\,\text{and}\,\,\pa$$ implies $i_1=k$ and $i_2=k+1$. Now, if
$\pi$ starts with $k,k+1,i_3,\dots,i_{k-1},1,i_k$ where $k+2\leq
i_3<\cdots<i_k\leq n$, then note that no occurrence of
$\tau\in\PP_k$ can contain the entries $k$ or $k+1$; further, an
occurrence of $321$ containing $k$ will exist in $\pi$ if and only
if there is such an occurrence containing $k+1$. Using this and
Lemma~\ref{l22},
\begin{eqnarray*}
f_k(n;k,k+1,i_3,\dots,i_{k-1},1,i_k)&=&f_k(n-1; k,i_3-1,\dots,i_{k-1}-1,1,i_k-1)\\
&=&f_k(n-2;k-1,i_3-2,\dots,i_{k-1}-2,i_k-2)\\
&=&f_k(n-3;i_3-3,\dots,i_{k-1}-3,i_k-3).\label{eqsa2}
\end{eqnarray*}
Combining the last equality with (\ref{eqsa}), (\ref{eqaa1}) and
the definitions of $C_d(n,m)$ and $A_d(n,m)$, we obtain the
desired equality
        $$C_k(n,k-1)=C_{k-1}(n-1,k-1)+A_{k-1}(n-3,k-2).$$

For (b), by definitions of $C_d(n,m)$ and $D_d(n,m)$, we have
\[C_{k-1}(n-1,k-1)=C_{k}(n-1,k-1)+D_{k-1}(n-1,k-2).\]
Combining with (a), it is enough to show
$D_{k-1}(n-1,k-2)=A_{k-2}(n-3,k-2)$. To this end, note that if
$\pi\in {\mathcal F}_k(n-1;k-1,i_1,\dots,i_{k-2},1)$ where
$k+1\leq i_1<\cdots<i_{k-2}\leq n-1$, then by Lemma~\ref{l22}
\[f_k(n-1;k-1,i_1,\dots,i_{k-2},1)=f_k(n-3;i_1-2,\dots,i_{k-2}-2).\]
By definitions of $D_d(n,m)$ and $A_d(n,m)$, we obtain the
required equality.
\end{proof}

\begin{proposition}\label{prob2}
Let $n\geq k-1$. The sequences $A_k, B_k, C_k$ and $D_k$ satisfy
the relations:
\begin{itemize}
\item[(a)] $A_k(n,k)=B_k(n-1,k-2)$;

\item[(b)] $B_k(n-1,k-2)=D_k(n-1,k-2)+B_k(n-1,k-1)$;

\item[(c)] $B_k(n-1,k-1)=B_k(n-2,k-2)$;

\item[(d)] $D_k(n-1,k-2)=A_{k-2}(n-4,k-2)+A_{k-1}(n-4,k-2)$;

\item[(e)] $A_k(n,k)-A_k(n-1,k)=A_{k-2}(n-4,k-2)+A_{k-1}(n-4,k-2)$.
\end{itemize}
\end{proposition}
\begin{proof}
For (a), if $\pi\in F_k(n;i_1,i_2,\dots,i_k)$ so that $k\leq
i_1<i_2<\cdots<i_k\leq n$, then avoidance of 321,
\[\pd\,\,\text{and}\,\,\pb\]
implies $i_1=k$ and $i_2=k+1$. Since no occurence of
$\tau\in\PP_k$ in $\pi$ can contain both entries $k$ and $k+1$, it
follows that
$f_k(n;k,k+1,i_3,\dots,i_k)=f_k(n-1;k,i_3-1,\dots,i_k-1)$, and
hence (a).

Part (b) is obtained similarly as in Lemma~\ref{l22}(b).

For (c), if $\pi\in  F_k(n-1;k,i_1,\dots,i_{k-1})$ where $k+1\leq
i_1<\cdots<i_{k-1}\leq n-1$, then avoidance of
$$321\,\,\text{and}\,\,\pb$$ implies again $i_1=k+1$. As in (a),
$f_k(n-1;k,k+1,i_2,\dots,i_{k-1})=f_k(n-2;k,i_2-1,\dots,i_{k-1}-1)$,
and (c) follows.

For (d), consider $D_k(n-1,k-2)$ along with a similar argument to
the one in Lemma~\ref{l22}(b):
\begin{equation*}
D_k(n-1,k-2)=\!\!\!\!\!\sum_{k+1\leq i_1<\cdots < i_{k-2}\leq n}
\!\!\!\!\!\!f_k(n;k,i_1,\dots,i_{k-2},1,2)+
\!\!\!\!\!\sum_{k+1\leq i_1<\cdots < i_{k-1}\leq n}
\!\!\!\!\!\!f_k(n;k,i_1,\dots,i_{k-2},1,i_{k-1}),
\end{equation*}
Part (d) follows from here as in the proof of (\ref{eqaa1}) and
(\ref{eqsa2}).

Finally, (a)--(d) yield (e).
\end{proof}

\subsection{Proof of Theorem~\ref{thm}}\label{third}
By (\ref{eqb1}), $A_k(n,k-1)-A_k(n,k)=C_k(n,k-1)$. Replacing $n$
by $n-1$ and subtracting yields
$$\big(A_k(n,k-1)-A_k(n,k)\big)-\big(A_k(n-1,k-1)-A_k(n-1,k)\big)=C_k(n,k-1)-C_k(n-1,k-1),$$
$$\Rightarrow\,\,
A_k(n,k-1)-A_k(n-1,k-1)=A_k(n,k)-A_k(n-1,k)+C_k(n,k-1)-C_k(n-1,k-1).$$
By Proposition \ref{prob1}(b) and Proposition \ref{prob2}(e),
\begin{eqnarray}
A_k(n,k-1)-A_k(n-1,k-1)&=&A_{k-2}(n-3,k-2)+A_{k-1}(n-3,k-2)\label{eqeq}\\
                     &+&A_{k-2}(n-4,k-2)+A_{k-1}(n-4,k-2).\nonumber
\end{eqnarray}
The result of Corollary~\ref{thss1} completes the proof of
Theorem~\ref{thm}.\qed

\subsection{Proof of Corollary~\ref{thmm}}
Since $f_k(n)=c_n=\frac{1}{n+1}\binom{2n}{n}$ for all
$n=0,1,\cdots,2k-1$, we have
$$\sum_{n\geq 0}
L_n^k(s)x^n=\sum_{i=0}^s(-x)^i\binom{s-i}{i}\left(f_k(x)-\sum_{j=0}^{s-i-1}x^jc_j\right).$$
Recall that the Chebyshev polynomials of the second kind satisfy
the relation
$$x^{s/2}U_s\ttx=\sum_{i=0}^s(-x)^i\binom{s-i}{i}$$
(see \cite[page~75-76]{Ri}), while the Catalan numbers satisfy the
relation
$$\sum_{i=0}^l(-1)^i\binom{s-i}{i}c_{l-i}=(-1)^l\binom{s-1-l}{l}$$
for all $l\leq  s-1$ (see \cite[page 152-154]{Ri}), and hence
$$\sl_{i=0}^s(-x)^i\binom{s-i}{i}\sum_{j=0}^{s-i-1}x^jc_j=\sl_{l=0}^{s-1}x^l\sum_{i=0}^l(-1)^i\binom{s-i}{i}c_{l-i}=\sum_{l=0}^{s-1}(-x)^l\binom{s-1-l}{l}.$$
Therefore
    $$\sum_{n\geq 0}L_n^k(s)x^n=x^{s/2}U_s\ttx-x^{(s-1)/2}U_{s-1}\ttx.$$
Finally, the Chebyshev polynomials of the second kind satisfy also
the relation
$$U_m\ttx=\frac{1}{\sqrt{x}}U_{m-1}\ttx-U_{m-2}\ttx$$
for all $m\geq 2$, hence by Theorem~\ref{thm} we get the desired
result. \qed
\section{Further results}
In this section we describe several directions which generalize
Theorem~\ref{thm} by use of the same arguments as in the proofs of
Theorem~\ref{thm} and Corollary~\ref{thmm}.

\subsection{First generalization}
For any $1\leq d \leq k-2$, let $X_{k,d}^1$ be the set of all
patterns
$$(d+1,d+2,\dots,k-1,1,2,\dots,j,k,j+1,j+2,\dots,d)$$
for $j=0,1,2,...,d$, plus the pattern $321$. For example,
$X_{3,1}^1=\{321, 213, 231\}$. Denote the number of permutations
in $ S_n(321,X_{k,d}^1)$ by $x_{k,d}^1(n)$.

\begin{theorem}\label{ext1}
Let $k\geq 3$ and $1\leq d\leq k-2$. For all $n\geq k$,
$$\sum_{j=0}^{k-1}(-1)^j\binom{k-1-j}{j}x_{k,d}^1(n-j)=\binom{n-d-1}{k-d-2},$$
and for all $0\leq n\leq k-1$,
$x_{k,d}^1(n)=\frac{1}{n+1}\binom{2n}{n}$. Thus, the generating
function for $\{x_{k,d}^1(n)\}_{n\geq 0}$ is
$$\sum_{n\geq 0}x_{k,d}^1(n)x^n=\dfrac{U_{k-2}\ttx+\frac{x^{k/2}}{(1-x)^{k-d-1}}}{\sqrt{x}U_{k-1}\ttx}\cdot$$
\end{theorem}
\begin{proof}
Similarly as before, we define $\mathcal{A}_s(n,m)$ to be the size
of $S_n(321,X_{k,d}^1)$ and obtain that
\begin{equation}
\mathcal{A}_{d+1}(n,k-d-1)=\sum_{j=0}^{k-1}(-1)^j\binom{k-1-j}{j}x_{k,d}^1(n-j)
\label{gen1}
\end{equation}
On the other hand, it is easy to verify that
$$\mathcal{A}_{d+1}(n,k-d-1)=\sum_{d+1\leq i_1<i_2<\cdots < i_{k-d-1}
\leq n}x_{k,d}^1(n;i_1,\dots,i_{k-d-1},1).$$ Since our
permutations avoid $X_{d,k}^1$, we set $i_{k-d-1}=n$. Thus, if
$\pi\in X_{k,d}^1(n;i_1,\dots,i_{k-d-2},n,1)$ with $d+1\leq
i_1<\cdots<i_{k-d-2}\leq n-1$, then $\pi_{k-d+1}<\cdots<\pi_n$.
This means
\begin{equation}
\mathcal{A}_{d+1}(n,k-d-1)=\binom{n-d-1}{k-d-2} \label{gen2}
\end{equation}
Combining (\ref{gen1}) and (\ref{gen2}) yields the desired
recursive relation. The rest of the theorem is easy to derive by
use of the same argument as in the proof of Corollary~\ref{thmm}.
\end{proof}

\begin{example}
{\rm  For $d=1$ and $k=3$, Theorem~\ref{ext1} yields}
\[x_{3,1}^1(n)-x_{3,1}^1(n-1)=1,\,\,{\rm and}\,\,\,
x_{3,1}^1(0)=x_{3,1}^1(1)=1,\,\,x_{3,1}^1(2)=2.\] {\rm Hence,
$x_{3,1}^1(n)=n$ for all $n\geq 1$ (see \cite{SS}). For $d=1$,
$k=4$ and $n\geq 0$, Theorem~\ref{ext1} yields}
    $$x_{4,1}^1(n)=|S_n(321,2341,2314)|=2^n-n,$$
{\rm while for $d=2$, $k=4$ and $n\geq 2$:}
    $$x_{4,2}^1(n)=|S_n(321,3412,3142,3124)|=3\cdot2^{n-2}-1.$$
\end{example}

\subsection{Second generalization} Let $X_{k,d}^2$ consist of the three
patterns $321$, $(d+1)(d+2)\cdots(k-1)1k23\cdots d$, and
$(d+1)(d+2)\cdots(k-1)k12\cdots d$. The number of
$X_{k,d}^2$-avoiding permutations in $ S_n$ is denoted by
$x_{k,d}^2(n)$. Similarly as in Theorem~\ref{ext1}, along with the
proof arguments of Theorem~\ref{thm} and Corollary~\ref{thmm}), we
arrive at
\begin{theorem}\label{ext2}
Let $k\geq 4$ and $2\leq d\leq k-2$. Then for all $n\geq  k$,
$$\begin{array}{l}
\sl_{j=0}^{k-1}(-1)^j\binom{k-1-j}{j}x_{k,d}^2(n-j)=
\sl_{j=0}^{k-1}(-1)^j\binom{k-1-j}{j}x_{k,d}^2(n-1-j)+\sl_{j=0}^{k-4}(-1)^j\binom{k-4-j}{j}x_{k,d}^2(n-3-j),
\end{array}$$
and $x_{k,d}^2(n)=\frac{1}{n+1}\binom{2n}{n}$ for all $0\leq n\leq
k-1$. Thus, the generating function for $\{x_{k,d}^2(n)\}_{n\geq
0}$ is
\[\sum_{n\geq 0}x_{k,d}^2(n)x^n=\dfrac{U_{k-1}\ttx+xU_{k-3}\ttx}{\sqrt{x}\left[U_k\ttx+xU_{k-2}\ttx\right]}\cdot\]
\end{theorem}

\begin{example}
{\rm For $d=2$, $k=4$ and $n\geq 3$, Theorem~\ref{ext2} yields}
\[x_{4,2}^2(n)=3x_{4,2}^2(n-1)-2x_{4,2}^2(n-2)+x_{4,2}^2(n-3),\,\,
{\rm and}\,\,x_{4,2}^2(0)=x_{4,2}^2(1)=1,\,\,x_{4,2}^2(2)=2.\]
\end{example}
A comparison of Theorem~\ref{ext2} for different cases suggests
that there should exist a bijection between the sets $
S_n(X_{k,2}^2)$ and $ S_n(132,X_{k,d}^2)$ for any $d$ such that
$2\leq d\leq k-2$. Producing such a bijection remains yet an open
question.

\subsection{Third generalization} Let $X_{k,d}^3$ consist of the two
patterns $321$ and $(d+1)(d+2)\cdots(k-1)1k23\cdots d$. The number
of $X_{k,d}^3$-avoiding permutations in $ S_n$ is denoted by
$x_{k,d}^3(n)$. Similarly as in the argument proofs of the main
theorem in \cite{MV3} and Theorem~\ref{thm}, we obtain

\begin{theorem}\label{ext3}
Let $k\geq  4$ and $2\leq  d\leq  k-2$. Then for all $n\geq  k$,
    $$\sum_{i=0}^k (-1)^j\binom{k-j}{j} x_{k,d}^2(n-j)=0,$$
and $x_{k,d}^2(n)=\frac{1}{n+1}\binom{2n}{n}$ for all $0\leq n\leq
k-1$. Thus, the generating function for $\{x_{k,d}^3(n)\}_{n\geq
0}$ is
            $$\sum_{n\geq 0}x_{k,d}^2(n)x^n=\dfrac{U_{k-1}(t)}{\sqrt{x}U_k\ttx}\cdot$$
\end{theorem}
Again, finding a direct bijection between the sets $
S_n(X_{k,d}^3)$ and $ S_n(132,(d+1)\cdots k12\cdots d)$ for any
$k\geq 4$ and $2\leq d\leq k-2$ is, as of now, an open question.


\begin{thebibliography}{WWW}
\bibitem[Bo1]{B1}
M.~B\'ona, The permutation classes equinumerous to the smooth
class, {\em Electron. J. Combin.} {\bf 5} (1998) \#R31.

\bibitem[Bo2]{B2}
M.~B\'ona, The solution of a conjecture of Stanley and Wilf for
all layered patterns, {\em J. Combin. Theory Ser. A}, {\bf 85}
(1999) 96--104.

\bibitem[BW]{BW}
S.~C.~Billey and G.~S.~Warrington, Kazhdan-lusztig polynomials for
$321$-hexagon-avoiding permutations,

\bibitem[CW]{CW}
T.~Chow and J.~West, Forbidden subsequences and Chebyshev
polynomials, {\em Discr. Math.} {\bf 204} (1999) 119--128.

\bibitem[Kn]{Kn}
D.E.~Knuth, The Art of Computer Programming, 2nd ed. Addison
Wesley, Reading, MA (1973).

\bibitem[Km]{Km}
D.~Kremer, Permutations with forbidden subsequences and a
generalized Schr\"oder number, {\em Discr. Math.} {\bf 218} (2000)
121--130.

\bibitem[Kr]{Kr}
C.~Krattenthaler, Permutations with restricted patterns and Dyck
paths, {\em Adv. in Applied Math.} {\bf 27} (2001), 510--530.

\bibitem[MV1]{MV1}
T.~Mansour and A.~Vainshtein, Restricted permutations, continued
fractions, and Chebyshev polynomials {\em Electron. J. Combin.}
\textbf{7} (2000) \#R17.

\bibitem[MV2]{MV2}
T.~Mansour and A.~Vainshtein, Restricted 132-avoiding
permutations, {\em Adv. Appl. Math.}  {\bf 126} (2001), 258--269.

\bibitem[MV3]{MV3}
T.~Mansour and A.~Vainshtein, Layered restrictions and Chebychev
polynomials (2000), {\em Annals of Combinatorics} {\bf 5} (2001),
451--458.

\bibitem[MV4]{MV4}
T.~Mansour and A.~Vainshtein, Restricted permutations and
Chebyshev polynomials, {\em S\'eminaire Lotharingien de
Combinatoire} {\bf 47} (2002), Article B47c.

\bibitem[R]{R}
A.~Robertson, Permutations containing and avoiding 123 and 132
patterns, {\em Discrete Mathematics and Theoretical Computer
Science}, {\bf 3} (1999) 151--154.

\bibitem[Ri]{Ri}
Th.~Rivlin, Chebyshev polynomials. From approximation theory to
algebra and number theory, John Wiley, New York (1990).

\bibitem[RWZ]{RWZ}
A.~Robertson, H.~Wilf, and D.~Zeilberger, Permutation patterns and
continuous fractions, {\em Electron. J. Combin.} {\bf 6} (1999)
\#R38.

\bibitem[SS]{SS}
R.~Simion, F.W.~Schmidt, Restricted Permutations, {\em Europ. J.
of Combinatorics} {\bf 6} (1985), 383--406.

\bibitem[SW]{SW}
Z.~Stankova and J.~West, Explicit enumeration of
$321$-hexagon-avoiding permutations, {\em Disc. Math.}, to appear.

\bibitem[W]{W}
J.~West, Generating trees and forbidden subsequences, {\em Discr.
Math.} {\bf 157} (1996), 363--372.
\end{thebibliography}
\end{document}